\documentclass[11pt]{article} 
\usepackage{asa}
\usepackage{harvard}
\usepackage{geometry}
\usepackage{fullpage}
\usepackage{epsfig}
\usepackage{mystyle}
\begin{document}
\newcommand{\mbf}[1]{\mbox{\boldmath $#1$}}
\def\al{\alpha}
\def\be{\beta}
\def\ga{\gamma}
\def\de{\delta}
\def\ep{\epsilon}
\def\ze{\zeta}
\def\et{\eta}
\def\th{\theta}
\def\ka{\kappa}
\def\la{\lambda}
\def\rh{\rho}
\def\si{\sigma}
\def\ta{\tau}
\def\up{\upsilon}
\def\ph{\phi}
\def\ch{\chi}
\def\ps{\psi}
\def\om{\omega}
\def\epv{\vareplison}
\def\thv{\vartheta}
\def\Ga{\Gamma}
\def\De{\Delta}
\def\Th{\Theta}
\def\La{\Lambda}
\def\Si{\Sigma}
\def\tSi{\tilde{\Sigma}}
\def\tS{\tilde{S}}
\def\tM{\tilde{M}}
\def\tOm{\tilde{\Omega}}
\def\Up{\Upsilon}
\def\Ph{\Phi}
\def\Ps{\Psi}
\def\Om{\Omega}
\def\del{\partial}
\def\beq{\begin{equation}}
\def\eeq{\end{equation}}
\def\mpar{\marginpar}
\def\lrarrow{\leftrightarrow}

\def\hk{{\mathcal H}_K}

\def\ss{\mbf{s}}
\def\cB{{\cal B}}
\def\cT{{\cal T}}
\def\cH{{\cal H}}
\def\cM{{\cal M}}
\def\cE{{\cal E}}
\def\cN{{\cal N}}
\def\cM{{\cal M}}
\def\cX{{\cal X}}
\def\cU{{\cal U}}
\def\cL{{\cal L}}
\def\cS{{\cal S}}
\def\cA{{\cal A}}
\def\cG{{\cal G}}
\def\cI{{\cal I}}

\newcommand{\bo}[1]{\mbox{{\boldmath{#1}}}}
\newcommand{\ms}{~~~}
\newtheorem{prop}{Proposition}
\newtheorem{Theorem}{Theorem}
\newtheorem{Lemma}{Lemma}

\newcommand{\ld}{\lambda}
\newcommand{\CH}{ {\cal H} }
\newcommand{\CM}{ {\cal M} }
\newcommand{\CG}{ {\cal G} }
\newcommand{\CL}{ {\cal L} }
\newcommand{\CN}{ {\cal N} }
\newcommand{\Pm}{ {\cal P}_{m-1} }
\newcommand{\CS}{ {\cal S} }
\newcommand{\CT}{ {\cal T} }

\newcommand{\bfa}{\mbox{\boldmath $a$}}
\newcommand{\bfbeta}{\mbox{\boldmath $\beta$}}
\newcommand{\bfc}{\mbox{\boldmath $c$}}
\newcommand{\bfd}{\mbox{\boldmath $d$}}
\newcommand{\bfe}{\mbox{\boldmath $e$}}
\newcommand{\bfdelta}{\mbox{\boldmath $\delta$}}
\newcommand{\bfeps}{\mbox{\boldmath $\epsilon$}}
\newcommand{\bff}{\mbox{\boldmath $f$}}
\newcommand{\bfF}{\mbox{\boldmath $F$}}
\newcommand{\bfrho}{\mbox{\boldmath $\rho$}}
\newcommand{\bfphi}{\mbox{\boldmath $\phi$}}
\newcommand{\bfs}{\mbox{\boldmath $s$}}
\newcommand{\bft}{\mbox{\boldmath $t$}}
\newcommand{\bftau}{\mbox{\boldmath $\tau$}}
\newcommand{\bfu}{\mbox{\boldmath $u$}}
\newcommand{\bfx}{\mbox{\boldmath $x$}}
\newcommand{\bfX}{\mbox{\boldmath $X$}}
\newcommand{\bfxi}{\mbox{\boldmath $\xi$}}
\newcommand{\bfy}{\mbox{\boldmath $y$}}
\newcommand{\bfz}{\mbox{\boldmath $z$}}
\newcommand{\bfmu}{\mbox{\boldmath $\mu$}}
\newcommand{\bfzero}{\mbox{\boldmath $0$}}
\newcommand{\bfone}{\mbox{\boldmath $1$}}

\newcommand{\hld}{\hat{\lambda}}
\newcommand{\hs}{\hat{\sigma}}
\newcommand{\fhtlht}{$\hat{f}_{\hat{\lambda}}$}
\newcommand{\ty}{\tilde{y}}
\newcommand{\bfty}{\mbox{\boldmath $\tilde{y}$}}
\newcommand{\teps}{\tilde{\epsilon}}
\newcommand{\bfteps}{\mbox{\boldmath $\tilde{\epsilon}$}}
\newcommand{\tf}{\tilde{f}}
\newcommand{\bftf}{\mbox{\boldmath $\tilde{f}$}}
\newcommand{\tty}{\tilde{\tilde{y}}}
\newcommand{\bftty}{\mbox{\boldmath $\tilde{\tilde{y}}$}}
\newcommand{\fl}{f_{\lambda, \theta}}
\newcommand{\gl}{g_{\lambda, \theta}}

\def\cite{\citeasnoun}
\thispagestyle{empty}
\begin{center}
{\large\bf
An Introduction to (Smoothing Spline) ANOVA Models in RKHS, With 
Examples in Geographical Data, Medicine, Atmospheric Science 
and Machine Learning}.
\end{center}
{\large\bf
\begin{center}
Grace Wahba \\
Department of Statistics \\
University of Wisconsin-Madison \\
{\tt http://www.stat.wisc.edu/\~{}wahba}
\end{center}
}
\begin{center}
This note has appeared in the Proceedings of the 13th 
IFAC Symposium on System Identification 2003, Rotterdam, 
549-559.
\end{center}
\section{Introduction}
Smoothing Spline ANOVA (SS-ANOVA) models in reproducing kernel
Hilbert spaces (RKHS) provide a very general framework for 
data analysis, modeling and learning in a variety 
of fields. Discrete, noisy scattered, direct and 
indirect observations can be accommodated with 
multiple inputs and multiple possibly correlated 
outputs and a variety 
of meaningful structures. 
The purpose of this paper is to 
give a brief overview of the approach and 
describe and contrast a series of applications, 
while noting some recent results.

\section{The general SS-ANOVA model}
The SS-ANOVA model
with Gaussian data has the form 
\begin{eqnarray}\label{gauss}
y_i = f(t_1(i),\cdots,t_d(i)) + \epsilon_i, ~~~~ i=1,\cdots,n,
\label{eq:obs-ssa}
\end{eqnarray}
where $\epsilon = (\epsilon_1, \cdots, \epsilon_n)'
 ~ \sim ~ N(0, ~ \sigma^2 I_{n \times n})$,
$t_{\al} \in \cT^{(\al)}$, where $\cT^{(\al)}$ is a measurable space,
$\al=1,\cdots,d;  (t_1,\cdots,t_d)= t 
\in  \cT = 
\cT^{(1)} \otimes \cdots \otimes \cT^{(d)}$, and 
$\si^2$ may be unknown. 
For $f$ satisfying 
some measurability conditions
a unique  ANOVA decomposition 
of $f$ of the form  
\begin{equation}\label{decomp}
f(t_1, \cdots , t_d) = \mu + 
\sum_{\al} f_{\al}(t_\al) + \sum_{\al\be}f_{\al\be}(t_{\al\be}) 
+ \cdots 
\end{equation}
can always be defined as follows:
Let $d\mu_{\al}$ be a probability measure on 
$\cT^{(\al)}$ and define the averaging operator 
$\cE_{\al}$ on $\cT$ by
\begin{equation}
({\cE}_{\al}f)(t)=
\int_{  {\cal{T}}^{(\al)}  }
f(t_1,\cdots,t_d)
d\mu_{\al}(t_\al).
\end{equation}
Then the identity is decomposed as
\beq\nonumber
I = \prod_{\al} (\cE_{\al} + (I-\cE_{\al})) =
\prod_{\al} \cE_{\al} + \sum_{\al}(I-\cE_{\al})\prod_{\be \neq \al} \cE_{\be} 
\eeq
\beq\label{E}
+ \sum_{\al < \be} (I-\cE_{\al})(I-\cE_{\be})
\prod_{\ga \neq \al,\be}\cE_{\ga}
+\cdots+ \prod_{\al}(I - \cE_{\al}).
\eeq
The components of this decomposition 
generate the ANOVA 
decomposition of $f$ of the form (\ref{decomp})
by  $C = (\prod_{\al}\cE_{\al})f, f_{\al} = ((I-\cE_{\al})\prod_{\be
\neq \al}\cE_{\be})f, f_{\al\be} = ((I -
\cE_{\al})(I-\cE_{\be})\prod_{\ga \neq \al,\be}\cE_{\ga})f $, 
and so forth. 
Further details
in the RKHS context may be found in 
\cite{wahba:1990}\cite{gu:wahba:1993}\cite{wahba:wang:gu:klein:klein:1995}

The idea behind SS-ANOVA is 
to construct an RKHS 
$\cH$  of functions on $\cT$ so that 
the components of the SS-ANOVA decomposition 
represent an orthogonal decomposition of $f$
in $\cH$. Then RKHS 
methods can be used to explicitly impose smoothness penalties 
of the form 
$\sum_{\al}\la_{\al}J_{\al}(f_{\al})
+ \sum_{\al\be}\la_{\al\be}J_{\al\be}(f_{\al\be}) + \cdots$, 
where, however, the series will be truncated at some point.
This is done as follows: 
Let $\CH^{(\al)}$ be an RKHS of functions on $\cT^{(\al)}$ with
$
\int _{\cT^{(\al)}} f_{\al} (t_{\al}) d\mu_{\al} = 0 \label{eq:side}
$
for $f_{\al} (t_{\al}) \in \CH^{(\al)}$,
and let  $[ 1^{(\al)}]$ be the one dimensional
space of constant functions on $\cT^{(\al)}$. 
Construct $\cH$ as 
\begin{equation}\nonumber
\cH  =  \prod _{j=1}^d (\{[ 1^{(\al)}] \} \oplus \{ \CH^{(\al)} \} )
\end{equation}
\begin{equation}
=  [1] \oplus \sum_j \CH^{(\al)} \oplus 
      \sum_{\al < \be} [\CH^{(\al)} \otimes \CH^{(\be)} ] 
      \oplus \cdots,
\label{eq:anova}
\end{equation}
where $[1]$ denotes the constant functions on $\cT$. 
With some abuse of 
notation, factors of the form 
$[1^{(\al)}]$ are omitted whenever they 
multiply a term of a different form. 
Thus  $\cH^{(\al)}$ is a shorthand for 
$[1^{(1)}]\otimes \cdots \otimes [1^{(\al-1)}]\otimes \cH^{(\al)}
\otimes [1^{(\al+1)}] \otimes \cdots \otimes[1^{(d)}]$ 
(which is a subspace 
of $\cH$). 
The components of the ANOVA decomposition are 
now in mutually orthogonal subspaces of $\cH$.
Note that the components will depend on the measures 
$d\mu_{\al}$ and these should be chosen in a specific 
application so that the fitted mean, main effects, two 
factor interactions, etc. have 
reasonable interpretations.  

Next,  $\cH^{(\al)}$ is decomposed into 
a parametric part and a smooth part,  by letting
$\cH^{(\al)} = \cH_{\pi}^{(\al)} \oplus \cH_{s}^{(\al)}$,
where $\cH_{\pi}^{(\al)}$ is finite dimensional (the ``parametric''
part) and $\cH_{s}^{(\al)}$
(the ``smooth'' part) is the orthocomplement of $\cH_{\pi}^{(\al)}$ in
$\cH^{(\al)}$.
Elements of $\cH_{\pi}^{(\al)}$ are not 
penalized through the device of letting 
$J_{\al}(f_{\al}) = \Vert P_{s}^{(\al)}f_{\al}\Vert^2 $
where $P_{s}^{(\al)}$ is the orthogonal projector 
onto $\cH_{s}^{(\al)}$.
$[\cH^{(\al)} \otimes \cH^{(\be)}]$ is now a direct sum of  four
orthogonal subspaces:
$[\cH^{(\al)} \otimes \cH^{(\be)}]=[\cH_{\pi}^{(\al)} \otimes \cH_{\pi}^{(\be)}]
\oplus [\cH_{\pi}^{(\al)} \otimes \cH_s^{(\be)}]
\oplus [\cH_s^{(\al)} \otimes \cH_\pi^{(\be)}]
\oplus [\cH_s^{(\al)} \otimes \cH_s^{(\be)}]$.
By convention the elements of the finite dimensional space
$[\cH_\pi^{(\al)}\otimes\cH_\pi^{(\be)}]$ will not be 
penalized.  
Continuing this way results in an orthogonal 
decomposition of $\cH$ into 
sums of
products of unpenalized finite dimensional subspaces, plus main
effects `smooth' subspaces, plus two factor interaction spaces of the form
parametric $\otimes$ smooth 
$[\cH_{\pi}^{(\al)} \otimes \cH_{s}^{(\be)}]$,
smooth $\otimes$ parametric
$[\cH_s^{(\al)} \otimes \cH_\pi^{(\be)}]$
and smooth $\otimes$ smooth
$[\cH_s^{(\al)} \otimes \cH_s^{(\be)}]$
and similarly for the
three and higher factor subspaces.

Now suppose that we have selected the model $\cM$, that is, we have
decided which subspaces will be included.  
Collect all of
the included unpenalized subspaces into a subspace, call it $\cH^0$,
of dimension $M$, and relabel the other subspaces as $\cH^{\be}, \be =
1,2,\cdots,p$.  
$\cH^{\be}$ may stand for a subspace
$\cH_{s}^{(\al)}$, or one of the three subspaces 
in the decomposition of 
$[\cH^{(\al)} \otimes \cH^{(\be)}]$ which contains 
at least one  `smooth' component, or, a higher order subspace with at 
least one `smooth' component. 
Collecting these subspaces as 
$\cM = \cH^0 \oplus \sum_{\be}\cH^{\be}$, 
the estimation problem in the Gaussian case becomes:
Find $f$ in $\cM={\cH}^0 \oplus\sum_{\be}\cH^{\be}$ to minimize
\beq\label{min}
\frac{1}{n} \sum_{i=1}^n (y_i - f(t(i)))^2 + 
\la \sum_{\be = 1}^p  \th_{\be}^{-1} \Vert P^{\be}f\Vert^2,
\eeq
where $P^{\be}$ is the orthogonal projector in $\cM$ onto $\cH^{\be}$,
and choose the (overparameterized) tuning parameters $\la, \th_{\be}$.
Bayesian confidence intervals, with the so-called `across the function'
property, are available for  these models. 

The residual sum of squares (RSS) in (\ref{min}) is replaced by 
the log likelihood 
\begin{equation}
{\cal{L}} (y,f) = -\sum_{i=1}^n[y_i f(t(i)) -b(f(t)))] 
\end{equation} 
for data from exponential families. Some of the examples 
below will involve Bernoulli $(0,1)$ data, in which case 
$b(f) = log(1+e^{f})$. Software for computing and tuning 
SS-ANOVA models may be found in the codes GRKPACK, 
RKPACK and gss and elsewhere, links to these and other 
spline related codes  
can be found via \\ 
{\tt http://www.stat.wisc.edu/\~{}wahba}
goto "SOFTWARE". Tuning methods are discussed 
in the first talk in this session. 
RSS may be replaced by 
robust functionals, or any convex functionals satisfying 
some mild conditions insuring uniqueness, 
and, in recent work on classification 
by support vector machines, RSS is replaced by so-called
hinge functions.

\section{Applications in Environmental Data}
\cite{gu:wahba:1993} considered data from 
the Eastern Lake Survey of 1984 which gave 
water acidity measurements and geographic 
locations, and other measurements
of lakes in the Blue Ridge Mountains
area. Of interest is the $pH$ as it depends on the 
geographic location and calcium concentration 
in the lakes. Model diagnostics were proposed
there, and the model
\begin{equation}
y_i = f_1(t_1(i)) + f_2(t_2(i)) + f_{1,2}(t_1(i), t_2(i)) + \ep_i
\end{equation}
was chosen, where $t_1$ is calcium content and $t_2$ is 
the pair (latitude, longitude). The thin plate spline 
penalty was imposed on the spatial variable. 
The calcium content and geography main effects 
models were plotted, and it can be seen that geography 
is a near proxy for elevation along the Blue 
Ridge mountains.

\section{Risk factor estimation}
\cite{wahba:wang:gu:klein:klein:1995} 
considered the risk of progression of diabetic 
retinopathy in a subpopulation of the 
Wisconsin Epidemiological Study of Diabetic
Retinopathy, whose baseline retinopathy 
score was below (i. e. good) a prespecified level.
The observations were 
$y_i = 1$ if the $i$th person's retinopathy 
progressed at the first followup, and 
$0$ if it had not. Here $f$ is the 
log odds ratio, $f = log[p/(1-p)]$. 
Three important variables were identified 
by informal means (see Section \ref{9})
and were $t_1 =$ duration of diabetes, 
$t_2 =$ glycosylated hemoglobin, and 
$t_3 =$ body mass index, and was modeled as 
\begin{equation}
f(t) = \mu + f_1(t_s) + a_2t_2 + 
f_3(t_3) + f_{13}(t_1, t_3).
\end{equation}
An interesting scientific result was found, 
that, persons in the study group with 
the longest duration of diabetes 
were at a lower risk, possibly because
they had survived longest without exceeding 
the prespecified threshold. 

\section{Time and Space Models on the Globe}
In \cite{wahba:luo:1997}\cite{luo:wahba:johnson:1997}
thirty years (1961-90) of Dec. Jan. Feb. average temperature 
measurements
at 1000 stations around the globe (with missing data)
was analyzed for spatial trends, as well as a global trend. 
Here 
$t = (t_1, t_2)= (x,P)$ where $x$ is year, and $P$ is 
(latitude, longitude).
The  RKHS of historical global temperature functions is
$\cH = [[1^{(1)}] \oplus [\phi]\oplus \cH_s^{(1)}] 
\otimes [[1^{(2)}] \oplus \cH_s^{(2)}]$, 
a collection of functions $f(x,P)$, on 
$\{1, 2,  ..., 30\} \otimes \cS$, where $\cS$ is 
the sphere, and $\cH$ and $f$ have 
corresponding decompositions given below:

\[
\begin{array}{rcccccccccccc}

\cH    &=&[1]     &\oplus&[\phi]  &\oplus&[\cH_s^{(1)}]
      &\oplus&[\cH_s^{(2)}]        &\oplus &[[\phi] \otimes \cH_s^{(2)}] 
      &\oplus &[\cH_s^{(1)} \otimes \cH_s^{(2)}]             \\

f(x,P) &=& C    &+     & d\ph(x)&+     & f_1(x)      &+     & f_2(P) &+&  \ph(x)f_{\ph,2}(P) &+& f_{12}(x,P) \\
\noalign{\vspace{-0pt}}
       &=& mean &+     &global  &+     & time        &+     &space   &+&   trend             &+& space-      \\       
\noalign{\vspace{-3pt}}
       & &      &      &time    &      & main        &      &main    & &   by~space          & & time         \\
\noalign{\vspace{-3pt}}
       & &      &     &trend    &      
       &\mbox{\it effect} & &\mbox{\it effect} & &\mbox{\it effect} & & interaction
\end{array}
\]
Here $\phi$ is a linear function which averages to $0$. A sum of squares 
of second differences was applied to the time variable, and a 
spline on the sphere penalty (\cite{wahba:1981d}\cite{wahba:1982a})
was applied to the space variable. For a cross country skiier in the 
midwest, as this author is, the results were very disappointing, 
in that they clearly showed a warming trend stretching from the 
midwest towards Alaska (trend by space term)  which was stronger than
the global mean trend. 
 
\section{Multiple correlated Bernoulli outcomes}
\cite{gao:wahba:klein:klein:2001} were  
motivated by a demographic
study involving a population with a variety 
of observed risk factors for several particular 
eye diseases, the outcomes were the incidence 
of one or more of several diseases or conditions 
in either or both of two eyes. Outcomes 
of the two eyes in a particular subject are 
presumed to be correlated, and incidences 
of the various outcomes may also be correlated. 
The amount of correlation may be of particular
interest. 
The risk factors could be person specific 
or eye-specific. The "two-eye" methods are 
a special case of what might be called  
"k-eye" methods where one person (unit) has several 
component outcomes which might have correlated 
outcomes, depending on unit-specific and component 
specific risk factors. 

The general log-linear model for multivariate Bernoulli data
goes as follows: 
Assuming there are $J$ different endpoints, and $K_j$ repeated
measurements for the $j$th endpoint, let $Y_{jk}$ denote the
$k$th measurement of the $j$th endpoint.
For example, in ophthalmological studies, we have
two repeated measurements for each disease: left eye and right eye.
In a typical longitudinal study, we have repeated measurements
over the time.
$Y=(Y_{jk},j=1,...,J,k=1,...,K_j)$ is
a multivariate Bernoulli outcome variable.
Let $X_{jk}=(X_{jk1},X_{jk2},...,X_{jkD})$
be a vector of predictor variables
ranging over the subset $\mathcal X$ of ${\mathcal R}^{D}$,
where $X_{jkd}$ denotes the $d$th predictor variable
for the $k$th measurement of the $j$th endpoint.
Some predictor variables may take different values
for different measurements
while others may be the same for all $Y_{jk}$'s.
For example, in ophthalmology studies, there may be present
both person-specific predictors and eye-specific predictors.
The person-specific predictors are the same for each person.
For the eye-specific predictors, the set of predictor variables is 
the same, but they may take different values for the left and right 
eyes. We can treat observations from both eyes as correlated repeated
measurements in our model.
Let $X=(X_{jk},j=1,...,J,k=1,...,K_j)$.
Then $(X,Y)$ is a pair of random vectors.
For a response vector $y=(y_{jk},j=1,...,J,k=1,...,K_j)$,
its joint probability distribution conditioning on
the predictor variables $X$ can be written as
\[
P(Y=y|X)=
\]
\[
\exp\{\sum_{j=1}^J\sum_{k=1}^{K_j} f_{jk} y_{jk}
       +\sum_{j=1}^J \sum_{k_1<k_2} \alpha_{jk_1,jk_2}  y_{jk_1} y_{jk_2}
\]
\[
+   \sum_{j_1<j_2} \sum_{k_1,k_2} \alpha_{j_1k_1,j_2k_2} 
                                   y_{j_1k_1} y_{j_2k_2}
       + ... + \alpha_{11,12,...,JK_J}  y_{11}y_{12}....y_{JK_J}
\]
\[
- b(f,\alpha)\},
\]
where
\begin{equation}\label{genmod}
b(f,\alpha)=
\end{equation}
\[
\log(1+\sum_{j,k}e^{f_{jk}}
  +\sum_{j_1,k_1}\sum_{j_2,k_2}
    e^{(f_{j_1k_1}+f_{j_2k_2}+\alpha_{j_1k_1,j_2k_2})}
  + ...
\]
\[
+ e^{(\sum_{all~ f} f + \sum_{all~ \alpha} \alpha)}).
\]
Let $M=\sum_{j=1}^{J} K_j$ be the length of the vector $Y$.
There are in total $2^M-1$ parameters:
$(f,\alpha)=(f_{11},f_{12},...,f_{JK_J},
             \alpha_{11,12},...,\alpha_{11,12,...,JK_J})$,
which may depend on $X$.
The parameter space is unconstrained. They have straightforward
interpretations in terms of conditional probabilities.
For example,
\begin{equation}
f_{jk}=logit(P(Y_{jk}=1 | Y^{(-jk)}=0,X))
\end{equation}
is the conditional logit function;
\begin{equation}
\alpha_{j_1k_1,j_2k_2}=\log OR(Y_{j_1k_1},Y_{j_2k_2} 
              | Y^{(-j_1k_1,-j_2k_2)}=0,X)
\end{equation}
is the conditional log odds ratio, which is a
meaningful way to measure pairwise association;
interpretations of other terms are given 
in the paper. 

$n$ independent observations
$(x_i,y_i),i=1,...,n$, are given, where $y_i=(y_{i11},y_{i12},...,y_{iJK_J})$ and 
$x_i=(x_{i11},x_{i12},...,x_{iJK_J})$. Here
$y_{ijk}$ and $x_{ijk}=(x_{ijk1},x_{ijk2},...,x_{ijkD})$
are the outcome variable and predictor vector for the
$k$th measurement of the $j$th endpoint of the $i$th
subject. 
Let $f_{jk}(i)$ be the conditional
logit function for the $k$th measurement of the $j$th endpoint
of the $i$th subject. 
There is little reason to believe the $f_{jk}$ will take different
functional forms for the same endpoint. Hence we can
assume $f_{ijk}=f_j(x_{ijk})$. The same reasoning applies to
the association terms. The $f_{jk}$ were modeled via 
SS-ANOVA in the paper, and a leaving-out-one-person 
based generalized cross validation for the smoothing 
parameters was obtained. 

\section{Multichotomous responses}\label{7}
\cite{linx:1998} considered multichotomous outcomes, the data is 
$(y_i,t(i))$ where $y_i$ is coded to show that the 
$i$ subject, with attribute vector $t(i)$ is in one of 
$k+1$ categories, $k > 1$. 
Let $p_j(t), j = 0, 1, \cdots, k$ be 
the probability that a subject with attribute vector 
$t$ is in category $k$, $\sum_{j=0}^k p_j(t) = 1.$
Let $f^j(t) = \log [p_j(t) /p_0(t)], j = 1, \cdots, k.$
Then 
\begin{eqnarray}
p_j(t) &=& \frac{e^{f^j(t)}}{1 + \sum_{j=1}^k e^{f^j(t)}}, j= 1, \cdots, k \\
p_0(t) &=& \frac{1}{1 + \sum_{j=1}^k e^{f^ j(t)}}.
\end{eqnarray}
The class label for the $i$th subject is coded as
$y_i = (y_{i1}, \cdots, y_{ik})$ where $y_{ij}=1$ if the 
$i$th subject is in class $j$ and $0$ otherwise. Letting 
$f = (f^1, \cdots, f^k)$ the negative log likelihood 
can be written as 
\begin{equation}\label{poly.lik}
\cL(y,f) = \sum_{i=1}^n \{ 
-\sum_{j=1}^k y_{ij}f^j(t_i) +
log(\sum_{j=1}^k 1+e^{f^j(t_i)})
\}.
\end{equation}
$f^j = \sum_{\nu_j=1}^M \phi_{\nu} + h^j$
where the $h^j$ 
can have an ANOVA decomposition.
Then $\la \Vert h \Vert_{\hk}^2$ in (\ref{min})
is replaced by \\
\begin{equation}\label{poly.J}
\sum_{j=1}^k
\sum_{\al}\la_{j\al} J_{j\al}( h^j_{\al}) +
\sum_{\al < \be} \la_{j\al\be}
J_{j\al\be}( h_{\al\be}^j) + \cdots.
\end{equation}
Ten year mortality 
data of a group of $n=646$ subjects 
with the risk factors 
age ($x_1$), glycosylated hemoglobin ($x_2$) 
and systolic blood pressure ($x_3$) were (among other things)
recorded at baseline and they were divided into 
four categories with respect to their status 
after ten years, as $0=$alive, $1=$ died of diabetes, 
$2=$died of heart disease, and $3=$died of other causes. 
Each of the $f^j , j = 1,2,3$ was modeled 
as $f^j(x_1, x_2, x_3) =
\mu^j + f_1^j(x_1) + f_2^j(x_2) + f_3^j(x_3) + f_{23}^j(x_2,x_3)$. 
The $p_j, j=0, \cdots, 3$ were estimated by minimizing
$\cI(y,f)=$ (\ref{poly.lik}) $+$ (\ref{poly.J}) 
and the multiple smoothing 
parameters estimated by a generalized cross validation 
method for polychotomous data given in 
\cite{linx:1998}. The plots graphically convey 
the suggestion that the younger deaths are disproportionately
diabetic, thus quickly raising further questions 
to confront the data base.

\section{The multicategory support vector machine}
The multicategory support vector machine (MSV) proposed 
in 
\cite{lee:lin:wahba:2002},\cite{lee:lin:wahba:2001b} 
considers the case where each subject is in one of 
$k$ categories labeled as $j = 1, \cdots, k$, 
as in the preceeding section, except for notational 
convenience there are $k$ instead of $k+1$ categories. 
The support vector machine is an efficient method 
for classification - it is not estimating 
the probability of membership in a particular 
category as before, but its target is an indicator 
as to which category as subject is in (or most likely
to be in)(see \cite{liny:2002}. 
The class label $y_i$ is now coded as a $k$ dimensional
vector with $1$ in the $j$th position if example $i$ 
is in category $j$ and $-\frac{1}{k-1}$ otherwise.
For example $y_i = (1, -\frac{1}{k-1}, \cdots, -\frac{1}{k-1})$
indicates that the $i$th example
is in category $1$. We define a $k$-tuple 
of separating functions
$f(t) = (f^1(t), \cdots f^k(t))$, with 
each $f^j = d^j + h^j$ with $h^j \in \hk$, 
and which will be 
required to satisfy 
a sum-to-zero constraint, $\sum_{j=1}^k f^j(t)= 0$, 
for all $t$ in $\cT$. Note that, unlike the estimate 
of Section \ref{7}, all categories are treated symmetrically.

Let $L_{jr} = 1, r \neq j$, $L_{jj}= 0, j, r = 1, \cdots, k$.  
Let $cat(y_i) =j$ if $y_i$ is from category $j$. 
Then, if $y_i$ is from category $j$,  
$L_{cat(y_i)r} = 0$ 
if $r = j$ and $1$ otherwise.
Then the 
MSVM is defined as the vector of functions
$f_{\la} = (f_{\la}^1, \cdots, f_{\la}^k)$, 
with each $h^k$ in $\hk$
satisfying the sum-to-zero constraint, 
which minimizes 
\begin{equation}\label{multicat}
\frac{1}{n}\sum_{i=1}^n \sum_{r=1}^k L_{cat(y_i)r}
(f^r(t_i) -y_{ir})_+ + \la\sum_{j=1}^k \Vert h^j\Vert_{\hk}^2.
\end{equation}
Generalizations of the penalty term are possible, if necessary.
It can be shown that the $k=2$ case reduces to the 
usual 2-category SVM just discussed, and
it is shown in \cite{lee:lin:wahba:2001b} that the target 
for the MSVM is 
$f(t) = (f^1(t), \cdots, f^k(t))$ with 
$f^j(t) = 1$ if $p_j(t)$ is bigger than the other $p_l(t)$ 
and $f^j(t) = -\frac{1}{k-1}$ otherwise. 
See also \cite{wahba:2002}.

\section{Summary}\label{9}
The SS-ANOVA models have proved to be useful in 
a  variety of modeling situations, 
only a few described here. In each case a 
tuning method which governs the bias-variance 
tradeoff must be employed, and, for very large 
sample sizes, efficient approximate methods need to be devised.
Model selection, that is, the determination of 
which variables and/or terms to include in the 
model is an important issue. 
\cite{zhang:wahba:lin:voelker:2001}\cite{zhang:wahba:lin:voelker:2002}
have recently
proposed likelihood basis 
pursuit, a nonparametric form of the LASSO, 
for the model 
selection problem associated with SS-ANOVA. 
Although a number of tuning methods for the 
various situations have been proposed, along 
with numerical methods for large data sets, 
a variety of problems remain to be investigated, 
including optimum nonlinear transformations 
of the variables, efficient computational 
methods, methods for covariates not missing 
at random, and public software for very large sample sizes
and for some of the more complex structures.

\bibliographystyle{agsm}
\bibliography{bib/nn/nn,bib/nn/nn1,bib/metbib/met,bib/metbib/met1,bib/metbib/met2,bib/gw11,bib/gw12,bib/gw1,bib/gw2,bib/gw3,bib/gw3b,bib/gw4,bib/gw5,bib/gw6,bib/gw7,bib/gw8,bib/cisbib/cis2,bib/k,bib/prep}

\begin{thebibliography}{xx}

\harvarditem[Gao et al.]{Gao, Wahba, Klein \&
  Klein}{2001}{gao:wahba:klein:klein:2001}
Gao, F., Wahba, G., Klein, R. \& Klein, B.  (2001), `Smoothing spline {ANOVA}
  for multivariate {B}ernoulli observations, with applications to ophthalmology
  data, with discussion', {\em J. Amer. Statist. Assoc.} {\bf 96},~127--160.

\harvarditem[Gu \& Wahba]{Gu \& Wahba}{1993}{gu:wahba:1993}
Gu, C. \& Wahba, G.  (1993), `Smoothing spline {ANOVA} with component-wise
  {B}ayesian ``confidence intervals''', {\em J. Computational and Graphical
  Statistics} {\bf 2},~97--117.

\harvarditem[Lee et al.]{Lee, Lin \& Wahba}{2001}{lee:lin:wahba:2001b}
Lee, Y., Lin, Y. \& Wahba, G.  (2001), Multicategory support vector machines,
  Technical Report 1043, Department of Statistics, University of Wisconsin,
  Madison WI.
\newblock To appear, {\em Computing Science and Statistics}, 33.

\harvarditem[Lee et al.]{Lee, Lin \& Wahba}{2002}{lee:lin:wahba:2002}
Lee, Y., Lin, Y. \& Wahba, G.  (2002), Multicategory support vector machines,
  theory, and application to the classification of microarray data and
  satellite radiance data, Technical Report 1063, Department of Statistics,
  University of Wisconsin, Madison WI.

\harvarditem[Lin]{Lin}{1998}{linx:1998}
Lin, X.  (1998), Smoothing spline analysis of variance for polychotomous
  response data, Technical Report 1003, PhD thesis, Department of Statistics,
  University of Wisconsin, Madison WI.
\newblock Available via G. Wahba's website.

\harvarditem[Lin]{Lin}{2002}{liny:2002}
Lin, Y.  (2002), `Support vector machines and the {B}ayes rule in
  classification', {\em Data Mining and Knowledge Discovery} {\bf 6},~259--275.

\harvarditem[Luo et al.]{Luo, Wahba \& Johnson}{1997}{luo:wahba:johnson:1997}
Luo, Z., Wahba, G. \& Johnson, D.  (1997), Spatial-temporal analysis of
  temperature using smoothing spline {ANOVA}, Technical Report 97-01,
  Pennsylvania State University Statistics Dept., State College PA.

\harvarditem[Wahba]{Wahba}{1981}{wahba:1981d}
Wahba, G.  (1981), `Spline interpolation and smoothing on the sphere', {\em
  SIAM J. Sci. Stat. Comput.} {\bf 2},~5--16.

\harvarditem[Wahba]{Wahba}{1982}{wahba:1982a}
Wahba, G.  (1982), `Erratum: Spline interpolation and smoothing on the sphere',
  {\em SIAM J. Sci. Stat. Comput.} {\bf 3},~385--386.

\harvarditem[Wahba]{Wahba}{1990}{wahba:1990}
Wahba, G.  (1990), {\em Spline Models for Observational Data}, SIAM.
\newblock CBMS-NSF Regional Conference Series in Applied Mathematics, v. 59.

\harvarditem[Wahba]{Wahba}{2002}{wahba:2002}
Wahba, G.  (2002), Soft and hard classification by reproducing kernel hilbert
  space methods, Technical Report 1067, Department of Statistics, University of
  Wisconsin, Madison WI.
\newblock to appear, Proceedings of the National Academy of Sciences.

\harvarditem[Wahba \& Luo]{Wahba \& Luo}{1997}{wahba:luo:1997}
Wahba, G. \& Luo, Z.  (1997), `Smoothing spline {ANOVA} fits for very large,
  nearly regular data sets, with application to historical global climate
  data', {\em Ann. Numer. Math.} {\bf 4},~579--597.

\harvarditem[Wahba et al.]{Wahba, Wang, Gu, Klein \&
  Klein}{1995}{wahba:wang:gu:klein:klein:1995}
Wahba, G., Wang, Y., Gu, C., Klein, R. \& Klein, B.  (1995), `Smoothing spline
  {ANOVA} for exponential families, with application to the {W}isconsin
  {E}pidemiological {S}tudy of {D}iabetic {R}etinopathy', {\em Ann. Statist.}
  {\bf 23},~1865--1895.
\newblock Neyman Lecture.

\harvarditem[Zhang et al.]{Zhang, Wahba, Lin, Voelker, Ferris, Klein \&
  Klein}{2001}{zhang:wahba:lin:voelker:2001}
Zhang, H., Wahba, G., Lin, Y., Voelker, M., Ferris, M., Klein, R. \& Klein, B.
  (2001), Variable selection via basis pursuit for non-{G}aussian data,
  Technical Report 1042, Statistics Department University of Wisconsin, Madison
  WI.
\newblock In Proceedings of the ASA Joint Statistical Meetings 2001 (CDROM),
  available from the American Statistical Association.

\harvarditem[Zhang et al.]{Zhang, Wahba, Lin, Voelker, Ferris, Klein \&
  Klein}{2002}{zhang:wahba:lin:voelker:2002}
Zhang, H., Wahba, G., Lin, Y., Voelker, M., Ferris, M., Klein, R. \& Klein, B.
  (2002), Variable selection and model building via likelihood basis pursuit,
  Technical Report 1059, Statistics Department University of Wisconsin, Madison
  WI.

\end{thebibliography}


\begin{thebibliography}{1}

\bibitem{wahba:1990}
G.~Wahba.
\newblock {\em Spline Models for Observational Data}.
\newblock SIAM, 1990.
\newblock CBMS-NSF Regional Conference Series in Applied Mathematics, v. 59.

\bibitem{gu:wahba:1993}
C.~Gu and G.~Wahba.
\newblock Smoothing spline {ANOVA} with component-wise {B}ayesian ``confidence
  intervals''.
\newblock {\em J. Computational and Graphical Statistics}, 2:97--117, 1993.

\bibitem{wahba:wang:gu:klein:klein:1995}
G.~Wahba, Y.~Wang, C.~Gu, R.~Klein, and B.~Klein.
\newblock Smoothing spline {ANOVA} for exponential families, with application
  to the {W}isconsin {E}pidemiological {S}tudy of {D}iabetic {R}etinopathy.
\newblock {\em Ann. Statist.}, 23:1865--1895, 1995.
\newblock Neyman Lecture.

\bibitem{luo:wahba:johnson:1998}
Z.~Luo, G.~Wahba, and D.~Johnson.
\newblock Spatial-temporal analysis of temperature using smoothing spline
  {ANOVA}.
\newblock {\em J. Climate}, 11:18--28, 1998.

\bibitem{lin:wahba:xiang:gao:klein:klein:2000}
X.~Lin, G.~Wahba, D.~Xiang, F.~Gao, R.~Klein, and B.~Klein.
\newblock Smoothing spline {ANOVA} models for large data sets with {B}ernoulli
  observations and the randomized {GACV}.
\newblock {\em Ann. Statist.}, 28:1570--1600, 2000.

\bibitem{gao:wahba:klein:klein:2001}
F.~Gao, G.~Wahba, R.~Klein, and B.~Klein.
\newblock Smoothing spline {ANOVA} for multivariate {B}ernoulli observations,
  with applications to ophthalmology data, with discussion.
\newblock {\em J. Amer. Statist. Assoc.}, 96:127--160, 2001.

\bibitem{wahba:2002}
G.~Wahba.
\newblock Soft and hard classification by reproducing kernel hilbert space
  methods.
\newblock Technical Report 1067, Department of Statistics, University of
  Wisconsin, Madison WI, 2002.
\newblock to appear, Proceedings of the National Academy of Sciences.

\bibitem{lee:lin:wahba:2002}
Y.~Lee, Y.~Lin, and G.~Wahba.
\newblock Multicategory support vector machines, theory, and application to the
  classification of microarray data and satellite radiance data.
\newblock Technical Report 1063, Department of Statistics, University of
  Wisconsin, Madison WI, 2002.

\bibitem{wahba:lin:leng:2002}
G.~Wahba, Y.~Lin, and C.~Leng.
\newblock Penalized log likelihood density estimation via smoothing-spline
  {ANOVA} and ran{GACV}-{C}omments to {H}ansen and {K}ooperberg, {S}pline
  adaptation in extended linear models.
\newblock {\em Statistical Science}, 17:33--37, 2002.

\end{thebibliography}

\end{document}